\documentclass[12pt]{article}
\usepackage{tikz}
\usetikzlibrary{tikzmark}
\usepackage{tikz-opm}
 \usepackage{amsmath, amsthm, amscd, amsfonts,graphicx}
\usepackage[matrix,arrow]{xy}

\newtheorem{thm}{Theorem}[section]
\newtheorem{dfn}{Definition}[section]
\newtheorem{lemma}{Lemma}[section]
\newtheorem{prop}{Proposition}[section]
\newtheorem{rem}{Remark}[section]
\newtheorem{ex}{Example}[section]
\newtheorem{prf}{Proof}[section]

\newcommand{\blem}{\begin{lemma}}
\newcommand{\elem}{\end{lemma}}

\newcommand{\bprf}{\begin{prf}}
\newcommand{\eprf}{\end{prf}}
\newcommand{\bthm}{\begin{thm}}
\newcommand{\ethm}{\end{thm}}

\newcommand{\bdf}{\begin{dfn}\rm}
\newcommand{\edf}{\end{dfn}}

\newcommand{\bpro}{\begin{prop}}
\newcommand{\epro}{\end{prop}}

\newcommand{\brem}{\begin{rem}\rm}
\newcommand{\erem}{\end{rem}}

\newcommand{\bex}{\begin{ex}\rm}
\newcommand{\eex}{\end{ex}}

\newcommand{\zerothemall}{
\setcounter{prop}{0} \setcounter{lemma}{0} \setcounter{rem}{0}
\setcounter{thm}{0} \setcounter{ex}{0} \setcounter{dfn}{0}
\setcounter{equation}{0} }

\title{A Model as a
Repeated Partnership Game
with  Discounting}
\author{ E. SOROURI, M. ESHAGHI GORDJI\\
Department of Mathematics, Semnan University,\\ P. O. Box 35195-363,
Semnan, Iran.\\
sorouri.e@semnan.ac.ir, meshaghi@semnan.ac.ir}
\date{}

\begin{document}
\maketitle

\begin{abstract}
\noindent
In this paper, we present a model of Partnership Game with respect to the important role of partnership and cooperation in nowdays life.
Since such interactions are repeated frequently, we study  this model as a Stage Game in the structure of infinitely repeated games with a discount factor $\delta$  and  Trigger strategy. We calculate and compare the payoffs of cooperation and violation and as an important result of this study, we show that each partner will adhere to the cooperation.

\vspace{5mm}

\noindent{\it  Keywords: game theory, Partnership Game, repeated game, Trigger strategy\\
JLE classification:C71, C73}

\end{abstract}

\section{Introduction}
\zerothemall
Since game theory examines situations in which decision-makers interact,  this theory has many applications such as firms competing for business, political candidates competing for votes, bidders competing in an auction, animals
fighting over prey, the arms race between countries, the relationship between parents and children, using the resources in nature,etc (see \cite{Bernheim},\cite{Bierman},\cite{Dixit},\cite{Gibbons},\cite{Osborne},\cite{OsborneM} and \cite{Petrosyan}).\\
On the other hand,  many of the strategic interactions in which we are involved are repeated interactions with the same people.  The relationship between the worker and the employer is  an example of this type. We can use the theory of repeated games to study   such behaviors. The main idea in this  theory is that a player may be deterred from exploiting her short-term advantage by the threat of punishment that reduces her long-term payoff.\\
In repeated  games with perfect  information that each player can observe the strategy used by other players, 
considering the discount factor $\delta$, it is possible that Nash equilibria of the repeated  game (supergame)  is  more efficient  than the Nash equilibria of the Stage Game, or one- period game. One of the important examples in this area is Cournots oligopoly game, which has been examined as an infinitely repeated game with  discount factor $\delta$ (see\cite{Abreu},\cite{Friedman},\cite{Rubinstein} and \cite{RubinsteinA}).\\

There are many activities and projects in which people contribute and the payoffs of those activities are derived from the efforts of each of the partners. Clearly, if any of the partners makes more efforts, more success will be achieved in these activities. But since more efforts by one person are beneficial to other people, they  may not have the motivation to work effectively on these projects. In fact, everyone chooses to make less effort and others to do more. With this view, another class of repeated games with imperfect  information has examined  models as Partnership Games with and without the discount factor.\cite{Fudenberg},\cite{Green},\cite{GreenE},\cite{Holmstrom},\cite{Porter},\cite{Radner}.

By getting the idea of Partnership Game in \cite{Lung} and \cite{Radner}, we have presented a more
complete model of participation and considering the role of collaboration in nowdays life and the fact that a collaborative activity can be repeated frequently, we study  the proposed model as a Stage Game in the
structure of infinitely repeated games with perfect  information and the discount factor $\delta$ between 0 and 1.
The results of this research encourage individuals to adhere to collaboration and cooperation, which is one of the most important goals of a social and modern life. 

\section{Model formulation and basic properties}
As a Complete Information Game, we assume that there is a collaborative activity
with two partners. The profit of this collaborative project depends on the
effort each partner spends on the project and is given by $\alpha(x _{1}+x_{2} +c _{1} (x_{1}x_{2}))$,
where $x_{1}$  is the amount of effort spent by partner 1 and $x_{2} $ is the amount
of effort spent by partner 2. Assume that $x_{1} ,x_{2} \in [0,\alpha] $. The value
$c_{1}\in [0,\frac{2}{\alpha}]$
 measures how complementary the efforts of the partners are. We
assume the amount of cost each player will incur for this effort is $c_{2}x_{i}^{2}$, where
$c_{2}\in [\frac{3}{2},2]$. Both players choose their effort independently and simultaneously,
and both want to maximize their share of the profit of the project which is
equally divided between two players. So the payoff function for partner i is

$$u_{i}(x_{1},x_{2})=\alpha(\frac{x_{1}+x_{2}}{2}+c_{1}(\frac{x_{1}x_{2}}{2}))-c_{2}x_{i}^{2}. $$

\section{Main results }
\subsection*{Nash equilibrium and the optimal amounts of effort}
Considering $\overline{x_{2}}$  as average effort, mathematical expectation, of the player 2 based on the belief of  player 1, we calculate the Nash equilibrium by finding the best response function  of each player
 \begin{equation}
 \dfrac{du_{1}(x_{1},\overline{x_{2}})}{dx_{1}}=\frac{\alpha}{2}+\frac{\alpha c_{1}}{2}\overline{x_{2}}-2c_{2}x_{1}=0,\label{3.2}
 \end{equation}
 
$$\dfrac{d^{2}u_{1}(x_{1},\overline{x_{2}})}{dx_{1}^{2}}=-2c_{2}<0.$$

Hence Equation  \ref{3.2} and  second derivative test  specify the best response function  of player1 as 
 $x_{1}=B_{1}(\overline{x_{2}})=\frac{\alpha}{4c_{2}}(1+c_{1}\overline{x_{2}})$.  
 Similarly the best response function  for player 2 is $x_{2}=B_{2}(\overline{x_{1}})=\frac{\alpha}{4c_{2}}(1+c_{1}\overline{x_{1}})$.
 A Nash equilibrium is a pair $(x_{1}^{*},x_{2}^{*})$ for which $x_{1}^{*}$ is a  best response to $x_{2}^{*}$ and 
$x_{2}^{*}$ is a  best response to $x_{1}^{*}$
\begin{equation}
\left\{ \begin{array}{cl}
x_{1}^{*}=B_{1}(x_{2}^{*})=\frac{\alpha}{4c_{2}}(1+c_{1}x_{2}^{*})\\
x_{2}^{*}=B_{2}(x_{1}^{*})=\frac{\alpha}{4c_{2}}(1+c_{1}x_{1}^{*}).
\end{array}\right.
\end{equation}
Solving these two equations, we find that $ x_{1}^{*}=x_{2}^{*}=\frac{\alpha}{4c_{2}-\alpha c_{1}}$.\\
The payoff of each player in the Nash equilibrium is 
\begin{equation}
 u_{i}(x_{1}^{*},x_{2}^{*})=\frac{\alpha^{2}}{2}(\frac{6c_{2}-\alpha c_{1}}{(4c_{2}-\alpha c_{1})^{2}}).
\end{equation}
We would like  to calculate the optimal amount of effort as follows
$$u(x_{1},x_{2})=\alpha(x_{1}+x_{2})+\alpha c_{1}(x_{1}x_{2})-c_{2}(x_{1}^{2}+x_{2}^{2}),$$
\begin{equation}
\left\{ \begin{array}{cl}
\dfrac{\partial u}{\partial x_{1}}=\alpha+\alpha c_{1}x_{2}-2c_{2}x_{1}=0\Longrightarrow x_{1}=\frac{\alpha(1+ c_{1}x_{2})}{2c_{2}}\\
\dfrac{\partial u}{\partial x_{2}}=\alpha+\alpha c_{1}x_{1}-2c_{2}x_{2}=0\Longrightarrow x_{2}=\frac{\alpha(1+ c_{1}x_{1})}{2c_{2}}.
\end{array}\right.
\end{equation}
By solving simultaneously the two equations $x_{1}=\frac{\alpha(1+ c_{1}x_{2})}{2c_{2}}$ and $x_{2}=\frac{\alpha(1+ c_{1}x_{1})}{2c_{2}}$, the result is $\widehat{x_{1}}=x_{1}=\dfrac{\alpha}{2c_{2}-\alpha c_{1}}$ and  $\widehat{x_{2}}=x_{2}=\dfrac{\alpha}{2c_{2}-\alpha c_{1}}$, 
clearly $x_{1}^{*}<\widehat{x_{1}}$ ,$x_{2}^{*}<\widehat{x_{2}}$.\\
On the other hand, according to  $c_{1}\in [0,\frac{2}{\alpha}]$  and $c_{2}\in [\frac{3}{2},2]$, we have 
\begin{equation}
D=\left\vert \begin{array}{cl}
-2c_{2} &\alpha c_{1}\\
\alpha c_{1} &-2c_{2}
\end{array}\right\vert=4c_{2}^{2}-\alpha^{2}c_{1}^{2}>0  \ ,\  \dfrac{\partial^{2} u}{\partial x_{1}^{2}}=-2c_{2}<0.
\end{equation}

So based on the second derivative test,  $(\widehat{x_{1}},\widehat{x_{2}})$ is a relative maximal point for $u(x_{1},x_{2})$.  With a simple calculation we have $u(\widehat{x_{1}},\widehat{x_{2}})=\frac{\alpha^{2}}{2c_{2}-\alpha c_{1}}$ . To determine the absolute maximum  of the optimal function, we must  compare values $u(0,0)=0$, $u(\alpha,\alpha)=\alpha^{2}(2-2c_{2}+\alpha c_{1})$ and $u(\widehat{x_{1}},\widehat{x_{2}})=\frac{\alpha^{2}}{2c_{2}-\alpha c_{1}}$. Clearly $u(0,0)<u(\alpha,\alpha)$ and 
$u(0,0)<u(\widehat{x_{1}},\widehat{x_{2}})$. Considering  $l=2c_{2}-\alpha c_{1}$ the relation
 $$\frac{\alpha^{2}}{2c_{2}-\alpha c_{1}}\geq\alpha^{2}(2-2c_{2}+\alpha c_{1})$$ is equivalent to$$\dfrac{1}{l}\geq (2-l)$$ that  is equivalent to$$(l-1)^{2}\geq 0,$$ which is always  true. So  in this case $u(x_{1},x_{2})$  has  the  absolute maximum in $$(\widehat{x_{1}},\widehat{x_{2}})=(\dfrac{\alpha}{2c_{2}-\alpha c_{1}},\dfrac{\alpha}{2c_{2}-\alpha c_{1}}).$$ Also 
\begin{equation}
        u_{i}(\widehat{x_{1}},\widehat{x_{2}})=\frac{\alpha^{2}}{2(2c_{2}-\alpha c_{1})}  \  \  for \ i=1,2.
\end{equation}

As an infinitely repeated game wiht Crime-Trigger strategy, we consider the Partnership model as a stage game in which  each of players has the same discount factor $\delta$.\\

$\mathbf{Theorem.1}$ If $\delta \in [\frac{ (4c_{2}-\alpha c_{1})^{2}}{8c_{2}(2c_{2}-\alpha c_{1})+(4c_{2}-\alpha c_{1})^{2}},1)$ then the Trigger strategy is  a Subgame Perfect Equilibruim, SPE .\\
Proof. We consider \begin{flushleft}
Stage game G:\\
Players:  \ Two players  \  $i=1,2 $ \\
Actions:  \ $ \forall i    \   x_{i}\in[0,\alpha]$\\
Stage Game Payoff: 
 $$ u_{i}(x_{1},x_{2})=\alpha(\frac{x_{1}+x_{2}}{2}+c_{1}(\frac{x_{1}x_{2}}{2}))-c_{2}x_{i}^{2}  \ for \ all \ players \ i. $$
 \end{flushleft}
We consider the Trigger strategy as follows

$$\  \forall i  \ S_{i}(h^{t})=\left\{ \begin{array}{rcl}
{\frac{\alpha}{2c_{2}-\alpha c_{1}}} &   \     \mbox{if}
& t=1 \\ \frac{\alpha}{2c_{2}-\alpha c_{1}} &   \   \mbox{if} & (\widehat{x_{1}},\widehat{x_{2}}),(\widehat{x_{1}},\widehat{x_{2}}),(\widehat{x_{1}},\widehat{x_{2}}),... \\
\frac{\alpha}{4c_{2}-\alpha c_{1}} & \   \mbox{if} & otherwise,
\end{array}\right.$$

where $h^{t}$ is the  history of game up to stage $t$. The concept of the above strategy is that in the first step, the amount of the effort of each player is  $\frac{\alpha}{2c_{2}-\alpha c_{1}}$.
If up to step $t-1$, each player has selected the amount $\frac{\alpha}{2c_{2}-\alpha c_{1}}$, then  the value $\frac{\alpha}{2c_{2}-\alpha c_{1}}$ is similarly chosen in step $t$, 
 otherwise the value of effort  of  the Nash equilibrium, $\frac{\alpha}{4c_{2}-\alpha c_{1}}$, will be selected.
We assume that the first  player adheres to the  above strategy. In order to determine the adherence of the second player to the above strategy , we will calculate her benefits from violations and non-violations.

First we presume that both players  adhere to the strategy. In this case the sequence of the players' selective combination will be as follows

$$(\frac{\alpha}{2c_{2}-\alpha c_{1}},\frac{\alpha}{2c_{2}-\alpha c_{1}}),(\frac{\alpha}{2c_{2}-\alpha c_{1}},\frac{\alpha}{2c_{2}-\alpha c_{1}}),... \ .$$

According to the above sequence, the payoff sequence for  the second player is
$$\frac{\alpha^{2}}{2(2c_{2}-\alpha c_{1})},\frac{\alpha^{2}}{2(2c_{2}-\alpha c_{1})},... \ .$$

Therefore the present value of the payoffs of the second  player is
$$\frac{\alpha^{2}}{2(2c_{2}-\alpha c_{1})}+\delta \frac{\alpha^{2}}{2(2c_{2}-\alpha c_{1})} +\delta^{2}\frac{\alpha^{2}}{2(2c_{2}-\alpha c_{1})}+...=\frac{\alpha^{2}}{2(2c_{2}-\alpha c_{1})}\frac{1}{1-\delta}.$$

Assuming that the first player  adheres to the strategy, we would like  to calculate the optimal amount of  the effort for  the second  player in case of the violation, $\mbox{max} \ u_{2}(\frac{\alpha}{2c_{2}-\alpha c_{1}},x_{2})$. In this case, we have
$$u_{2}(\frac{\alpha}{2c_{2}-\alpha c_{1}},x_{2})=\frac{\alpha}{2}(\frac{\alpha}{2c_{2}-\alpha c_{1}}+x_{2})+\frac{\alpha c_{1}}{2}(\frac{\alpha x_{2}}{2c_{2}-\alpha c_{1}})-c_{2}x_{2}^{2}$$
$$\dfrac{du_{2}}{dx_{2}}=\frac{\alpha}{2}+\frac{c_{1}\alpha^{2}}{2(2c_{2}-\alpha c_{1})}-2c_{2}x_{2}=0 \Longrightarrow
x_{2}=\frac{\alpha}{2(2c_{2}-\alpha c_{1})}$$
$$, \dfrac{d^{2}u_{2}}{dx_{2}}=-2c_{2}<0.$$

Therefore $x_{2}=\frac{\alpha}{2(2c_{2}-\alpha c_{1})}$ is a relative maximum for $u_{2}(\frac{\alpha}{2c_{2}-\alpha c_{1}},x_{2})$.\\
With a simple comparison between $u_{2}(\frac{\alpha}{2c_{2}-\alpha c_{1}},0)=\frac{\alpha^{2}}{2(2c_{2}-\alpha c_{1})}$,
$$u_{2}(\frac{\alpha}{2c_{2}-\alpha c_{1}},\alpha)=\frac{\alpha^{2}}{2(2c_{2}-\alpha c_{1})}+\frac{c_{2}\alpha^{2}}{(2c_{2}-\alpha c_{1})}(1-(2c_{2}-\alpha c_{1}))$$
and$$u_{2}(\frac{\alpha}{2c_{2}-\alpha c_{1}},\frac{\alpha}{2(2c_{2}-\alpha c_{1})})=(\frac{\alpha^{2}(5c_{2}-2\alpha c_{1})}{4(2c_{2}-\alpha c_{1})^{2}})=\frac{\alpha^{2}}{2(2c_{2}-\alpha c_{1})}+
  \frac{c_{2}\alpha^{2}}{4(2c_{2}-\alpha c_{1})^{2}}$$

 and  considering
   $$2c_{2}-\alpha c_{1}\geq 1,$$ it follows  that $x_{2}=\frac{\alpha}{2(2c_{2}-\alpha c_{1})}$ is an absolute maximal for $u_{2}(\frac{\alpha}{2c_{2}-\alpha c_{1}},x_{2})$.\\
   
The important point is
 $$   u_{2}(\widehat{x_{1}},\frac{\alpha}{2(2c_{2}-\alpha c_{1})})=u_{2}(\frac{\alpha}{2c_{2}-\alpha c_{1}},\frac{\alpha}{2(2c_{2}-\alpha c_{1})})=u_{2}(\widehat{x_{1}},\widehat{x_{2}})+
  \frac{c_{2}\alpha^{2}}{4(2c_{2}-\alpha c_{1})^{2}}$$ so  
   $u_{2}(\widehat{x_{1}},\frac{\alpha}{2(2c_{2}-\alpha c_{1})})>u_{2}(\widehat{x_{1}},\widehat{x_{2}})$  while $\frac{\alpha}{2(2c_{2}-\alpha c_{1})}<\widehat{x_{2}}$.\\
 This means that the second player can achieve more payoff with an effort less than the optimal amount of effort.  
   
   Let's consider the selection sequence of the players in case of the second player's violation as follows
$$(\frac{\alpha}{2c_{2}-\alpha c_{1}},\frac{\alpha}{2(2c_{2}-\alpha c_{1})}),(\frac{\alpha}{4c_{2}-\alpha c_{1}},\frac{\alpha}{4c_{2}-\alpha c_{1}}),(\frac{\alpha}{4c_{2}-\alpha c_{1}},\frac{\alpha}{4c_{2}-\alpha c_{1}}),... \ .$$

So  the sequence of the payoffs of   player 2 is as follows
$$\frac{\alpha^{2}}{2(2c_{2}-\alpha c_{1})}+
  \frac{c_{2}\alpha^{2}}{4(2c_{2}-\alpha c_{1})^{2}} \ , \ \frac{\alpha^{2}}{2}(\frac{6c_{2}-\alpha c_{1}}{(4c_{2}-\alpha c_{1})^{2}}) \ , \ 
  \frac{\alpha^{2}}{2}(\frac{6c_{2}-\alpha c_{1}}{(4c_{2}-\alpha c_{1})^{2}}) , ... \ .$$

Therefore in  the violation, the present value of the payoffs of the second  player is
$$(\frac{\alpha^{2}}{2(2c_{2}-\alpha c_{1})}+
  \frac{c_{2}\alpha^{2}}{4(2c_{2}-\alpha c_{1})^{2}}) +\delta (\frac{\alpha^{2}}{2}(\frac{6c_{2}-\alpha c_{1}}{(4c_{2}-\alpha c_{1})^{2}}) )+ \delta^{2}(\frac{\alpha^{2}}{2}(\frac{6c_{2}-\alpha c_{1}}{(4c_{2}-\alpha c_{1})^{2}})) + ...$$
$$=(\frac{\alpha^{2}(5c_{2}-2\alpha c_{1})}{4(2c_{2}-\alpha c_{1})^{2}})+\delta (\frac{\alpha^{2}(6c_{2}-\alpha c_{1})}{2(4c_{2}-\alpha c_{1})^{2}(1-\delta)}). $$

  Since player 2 will not violate if her  payoff is greater than or   at least equal to the non-violation, then we have to have
$$\frac{\alpha^{2}}{2(2c_{2}-\alpha c_{1})}\frac{1}{(1-\delta)}\geq(\frac{\alpha^{2}(5c_{2}-2\alpha c_{1})}{4(2c_{2}-\alpha c_{1})^{2}})+\delta (\frac{\alpha^{2}(6c_{2}-\alpha c_{1})}{2(4c_{2}-\alpha c_{1})^{2}(1-\delta)}).$$

 It is easy to check  that the above inequation is equivalent to 
 
 $$\delta\geq\frac{(4c_{2}-\alpha c_{1})^{2}}{(4c_{2}-\alpha c_{1})^{2}+8c_{2}(2c_{2}-\alpha c_{1})}.$$

 With the same process for player 1, one can show  that if  $\delta \in [\frac{ (4c_{2}-\alpha c_{1})^{2}}{8c_{2}(2c_{2}-\alpha c_{1})+(4c_{2}-\alpha c_{1})^{2}},1)$ then according to the  Trigger strategy, the players will continue the cooperation.  
 Therefore  the Trigger strategy is a SPE, that is, in each subgame the players choose the  cooperation, and no player intends to violate because his payoffs reduce  in comparison  with  cooperation. This completes the proof.

 $\mathbf{Theorem .2}$   
  In the partnership game for all $\delta \in (0,\frac{ (4c_{2}-\alpha c_{1})^{2}}{8c_{2}(2c_{2}-\alpha c_{1})+(4c_{2}-\alpha c_{1})^{2}})$  we can define the Trigger strategy, in which each player's level of effort, $\overline{x}$, is   greater than $x^{*}$ and less than $\widehat{x}$.\\
  Proof.We consider \begin{flushleft}
Stage game G:\\
Players:  \ Two players  \  $i=1,2 $ \\
Actions:  \ $ \forall i    \   x_{i}\in[0,\alpha]$\\
Stage Game Payoff:  
$$ u_{i}(x_{1},x_{2})=\alpha(\frac{x_{1}+x_{2}}{2}+c_{1}(\frac{x_{1}x_{2}}{2}))-c_{2}x_{i}^{2}  \ for \ all \ player \ i. $$
 \end{flushleft}
Also, we consider the Trigger strategy  as follows

$$   \  \forall i  \ S_{i}(h^{t})=\left\{ \begin{array}{rcl}
\overline{x} &   \     \mbox{if}
& t=1 \\ \overline{x} &   \   \mbox{if} & (\overline{x},\overline{x}),(\overline{x},\overline{x}),(\overline{x},\overline{x}),... \\
\frac{\alpha}{4c_{2}-\alpha c_{1}} & \   \mbox{if} & otherwise,
\end{array}\right.$$
where $h^{t}$ is the  history of game up stage $t$.
The concept of the above strategy is that in the first step, the amount of the effort of each player is  $\overline{x}$.
If up to step $t-1$, each player has selected the amount $\overline{x}$, then  the value $\overline{x}$ is similarly chosen in step $t$ 
 otherwise the value of effort  of  the Nash equilibrium, $\frac{\alpha}{4c_{2}-\alpha c_{1}}$, will be selected.\\
 
First we presume that both players  adhere to the strategy in which case the sequence of the players' selective combination will be as follows
$$(\overline{x} ,\overline{x}) , \  (\overline{x} ,\overline{x}) , \ (\overline{x} ,\overline{x}) ,...\ .$$

According to the above sequence, the payoff sequence for player 2  is 
$$u_{2}(\overline{x} ,\overline{x})  , \ u_{2}(\overline{x} ,\overline{x}) , \ u_{2}(\overline{x} ,\overline{x}) ,... \ ,$$

where  $u_{2}(\overline{x} ,\overline{x})=\alpha \overline{x}+ {\overline{x}}^{2}(\frac{\alpha c_{1}}{2}-c_{2})$.\\
Therefore the present value of the payoffs of   player 2  in  case of non-violation  is
$$(\alpha \overline{x}+ {\overline{x}}^{2}(\frac{\alpha c_{1}}{2}-c_{2}))+\delta (\alpha \overline{x}+ {\overline{x}}^{2}(\frac{\alpha c_{1}}{2}-c_{2})) +\delta^{2}(\alpha \overline{x}+ {\overline{x}}^{2}(\frac{\alpha c_{1}}{2}-c_{2})) +...$$
$$=(\alpha \overline{x}+ {\overline{x}}^{2}(\frac{\alpha c_{1}}{2}-c_{2}))\frac{1}{1-\delta}.$$

Assuming that player 1 selects the level of effort $\overline{x}$ and player 2 intends to violate  from $\overline{x}$, we calculate the optimal amount of effort, $x_{*}$, that maximizes  her payoff
$$u_{2}(\overline{x} ,x_{*})=\frac{\alpha}{2}(\overline{x} +x_{*})+\frac{\alpha c_{1}}{2}(\overline{x} x_{*})-c_{2}x_{*}^{2}$$
$$\dfrac{du_{2}}{dx_{*}} =\frac{\alpha}{2}+\frac{\alpha c_{1}}{2}\overline{x}-2c_{2} x_{*}=0  \Longrightarrow  x_{*}=\frac{\alpha}{4c_{2}}(1+c_{1}\overline{x})$$
$$,\dfrac{d^{2}u_{2}}{dx_{*}^{2}}=-2c_{2}<0.$$

Therefore, according to the second derivative test,  $ x_{*}=\frac{\alpha}{4c_{2}}(1+c_{1}\overline{x})$  is a relative maximal point for $u_{2}(\overline{x} ,x_{*})$.

 Also since  $c_{1}\in [0,\frac{2}{\alpha}] $ and $c_{2}\in [\frac{3}{2},2]$ it is easy to get
 $ x_{*}=\frac{\alpha}{4c_{2}}(1+c_{1}\overline{x})\leq\frac{\alpha}{2}$.
To determine the absolute maximum of $u_{2}(\overline{x} ,x_{*})$, we need to compare the values 
$u_{2}(\overline{x} ,0)=\frac{\alpha}{2}\overline{x}$,  $ u_{2}(\overline{x} ,\alpha)=\frac{\alpha}{2}(\overline{x}+\alpha (1+c_{1}\overline{x})-2c_{2}\alpha)$ and 
 $$u_{2}(\overline{x} ,x_{*})=\frac{\alpha}{2}(\overline{x}+\frac{\alpha}{8c_{2}}(1+c_{1}\overline{x})^{2}).$$
Clearly always $u_{2}(\overline{x} ,x_{*})>u_{2}(\overline{x} ,0)$.\\
On the other hand, 

$$\begin{array}{l}
u_{2}(\overline{x} ,x_{*})\geq u_{2}(\overline{x} ,\alpha) \\
\Longleftrightarrow\frac{\alpha}{2}(\overline{x}+\frac{\alpha}{8c_{2}}(1+c_{1}\overline{x})^{2})\geq \frac{\alpha}{2}(\overline{x}+\alpha (1+c_{1}\overline{x})-2c_{2}\alpha)\\
\Longleftrightarrow \frac{(1+c_{1}\overline{x})}{8c_{2}}+\frac{2c_{2}}{(1+c_{1}\overline{x})} \geq 1 \\
\Longleftrightarrow (4c_{2}-(1+c_{1}\overline{x}))^{2}\geq 0.
\end{array}$$

Because  $4c_{2}-(1+c_{1}\overline{x}))^{2}\geq 0$ is always true, then always
 $$u_{2}(\overline{x} ,x_{*})\geq u_{2}(\overline{x} ,\alpha).$$ 
Since players are always looking for less effort and more payoff, even if $u_{2}(\overline{x} ,x_{*})= u_{2}(\overline{x} ,\alpha)$ then player 2 always chooses less effort, $x_{*}=\frac{\alpha}{4c_{2}}(1+c_{1}\overline{x})$.\\
In the above argument, $x_{*}$  is  the optimal amount of effort for player 2 when player 1 selects  $\overline{x}$.\\
This way, the selection sequence of the players in case of the second player's violation is
$$(\overline{x} ,\frac{\alpha}{4c_{2}}(1+c_{1}\overline{x})) \ , \  (\frac{\alpha}{4c_{2}-\alpha c_{1}},\frac{\alpha}{4c_{2}-\alpha c_{1}}) \ , \ (\frac{\alpha}{4c_{2}-\alpha c_{1}},\frac{\alpha}{4c_{2}-\alpha c_{1}}) \ , \ ... \ .$$

According to the above sequence,  the sequence of the payoffs of   player 2 is as follows
$$\frac{\alpha}{2}(\overline{x}+\frac{\alpha}{8c_{2}}(1+c_{1}\overline{x})^{2}) \ , \ \frac{\alpha^{2}}{2}(\frac{6c_{2}-\alpha c_{1}}{(4c_{2}-\alpha c_{1})^{2}}) \ , \
 \frac{\alpha^{2}}{2}(\frac{6c_{2}-\alpha c_{1}}{(4c_{2}-\alpha c_{1})^{2}}) \ , \  ... \ . $$

So in case of violation, the present value of the payoffs of the second  player is
$$(\frac{\alpha}{2}(\overline{x}+\frac{\alpha}{8c_{2}}(1+c_{1}\overline{x})^{2}))+\delta \frac{\alpha^{2}}{2}(\frac{6c_{2}-\alpha c_{1}}{(4c_{2}-\alpha c_{1})^{2}}) + \delta^{2}\frac{\alpha^{2}}{2}(\frac{6c_{2}-\alpha c_{1}}{(4c_{2}-\alpha c_{1})^{2}})+
...$$
$$=(\frac{\alpha}{2}(\overline{x}+\frac{\alpha}{8c_{2}}(1+c_{1}\overline{x})^{2}))+\frac{\alpha^{2}}{2}(\frac{6c_{2}-\alpha c_{1}}{(4c_{2}-\alpha c_{1})^{2}})(\frac{\delta}{1-\delta}).$$

Obviously, player 2 will adhere to Trigger strategy if
$$(\alpha \overline{x}+ {\overline{x}}^{2}(\frac{\alpha c_{1}}{2}-c_{2}))\frac{1}{1-\delta}\geq
(\frac{\alpha}{2}(\overline{x}+\frac{\alpha}{8c_{2}}(1+c_{1}\overline{x})^{2}))+\frac{\alpha^{2}}{2}(\frac{6c_{2}-\alpha c_{1}}{(4c_{2}-\alpha c_{1})^{2}})(\frac{\delta}{1-\delta}).$$

By calculating, it is determined that the above inequality is equivalent to 
 $A \overline{x}^{2}+B\overline{x}+c\geq0$  in which 
 $$A=\frac{-1}{16c_{2}}((4c_{2}-\alpha c_{1})^{2}-\alpha^{2}c_{1}^{2}\delta)=-(8c_{2}(2c_{2}-\alpha c_{1})+\alpha^{2}c_{1}^{2}(1-\delta))<0,$$
 $$B=\frac{\alpha}{8c_{2}}((4c_{2}-\alpha c_{1})+\delta(4c_{2}+\alpha c_{1}))>0$$
 $$C=-\frac{\alpha^{2}}{16c_{2}}(\delta (\frac{32c_{2}^{2}-\alpha^{2}c_{1}^{2}}{(4c_{2}-\alpha c_{1})^{2}})+1)<0.$$

Put $p(\overline{x})=A \overline{x}^{2}+B\overline{x}+c$, then this equation has $\sqrt{\Delta} =\frac{2\alpha c_{2}\delta}{4c_{2}-\alpha c_{1}}$ and its roots are
$$\overline{x}_{1}=\frac{\alpha}{4c_{2}-\alpha c_{1}}   \ , \  
\overline{x}_{2}= \frac{\alpha}{4c_{2}-\alpha c_{1}}\frac{(4c_{2}-\alpha c_{1})^{2}-\delta \alpha^{2}c_{1}^{2}+32 \delta c_{2}^{2}}{(4c_{2}-\alpha c_{1})^{2}-\delta \alpha^{2}c_{1}^{2}}.$$

In which  $\overline{x}_{1}<\overline{x}_{2}$  and $\overline{x}_{1}$ is the Nash equilibrium. Also,

$$\begin{array}{l}
\delta <\frac{ (4c_{2}-\alpha c_{1})^{2}}{32c_{2}^{2}-16\alpha c_{1}c_{2}+\alpha^{2}c_{1}^{2}}\\
\Longleftrightarrow \frac{32\delta c_{2}^{2}}{(4c_{2}-\alpha c_{1})^{2}-\delta \alpha^{2}c_{1}^{2}}<\frac{2c_{2}}{2c_{2}-\alpha c_{1}}\\
\Longleftrightarrow   1+( \frac{32\delta c_{2}^{2}}{(4c_{2}-\alpha c_{1})^{2}-\delta \alpha^{2}c_{1}^{2}})<1+(\frac{2c_{2}}{2c_{2}-\alpha c_{1}})=\frac{4c_{2}-\alpha c_{1}}{2c_{2}-\alpha c_{1}}\\
\Longleftrightarrow \frac{\alpha}{4c_{2}-\alpha c_{1}}\frac{(4c_{2}-\alpha c_{1})^{2}-\delta \alpha^{2}c_{1}^{2}+32 \delta c_{2}^{2}}{(4c_{2}-\alpha c_{1})^{2}-\delta \alpha^{2}c_{1}^{2}}<\frac{\alpha}{2c_{2}-\alpha c_{1}}
\end{array}$$

So, according to the above calculations $\frac{\alpha}{4c_{2}-\alpha c_{1}}<\overline{x}_{2}<\frac{\alpha}{2c_{2}-\alpha c_{1}}$.\\
 By specifying the sign $p(\overline{x})$, it follows that $p(\overline{x})$ is always nonnegative between  two roots. On the other hand, if $x_{1} =x_{2} =\overline{x}$  then the calculations indicate that $u_{i}(\overline{x},\overline{x})=\alpha \overline{x} +(\overline{x})^{2}(\frac{\alpha c_{1}}{2}-c_{2})$ has two roots  of  $0$ and $\frac{2\alpha}{2c_{2}-\alpha c_{1}}$  and it is maximal in $\frac{\alpha}{2c_{2}-\alpha c_{1}}$. So, if $\frac{\alpha}{4c_{2}-\alpha c_{1}} <\overline{x}<\frac{\alpha}{2c_{2}-\alpha c_{1}}$  then  $u_{i}(\overline{x},\overline{x})$ is an increasing function.
So, if players choose the level of more effort, their payoffs will be greater. 
Therefore,  the purpose of solving $p(\overline{x})\geq0$ is the largest $\overline{x}$ for which  $p(\overline{x})\geq0$. Hence the highest value of  $\overline{x}$ is $\overline{x}_{2}$.\\

$\mathbf{Corollary.1}$. We have  $ \lim_{\delta\rightarrow 0}\overline{x}_{2}= \lim_{\delta\rightarrow 0}\frac{\alpha}{4c_{2}-\alpha c_{1}}\frac{(4c_{2}-\alpha c_{1})^{2}-\delta \alpha^{2}c_{1}^{2}+32 \delta c_{2}^{2}}{(4c_{2}-\alpha c_{1})^{2}-\delta \alpha^{2}c_{1}^{2}}
=\frac{\alpha}{4c_{2}-\alpha c_{1}}$\\

 $ ,\lim_{\delta\rightarrow \frac{ (4c_{2}-\alpha c_{1})^{2}}{8c_{2}(2c_{2}-\alpha c_{1})+(4c_{2}-\alpha c_{1})^{2}}}\overline{x}_{2}= \lim_{\delta\rightarrow \frac{ (4c_{2}-\alpha c_{1})^{2}}{8c_{2}(2c_{2}-\alpha c_{1})+(4c_{2}-\alpha c_{1})^{2}}}\frac{\alpha}{4c_{2}-\alpha c_{1}}\frac{(4c_{2}-\alpha c_{1})^{2}-\delta \alpha^{2}c_{1}^{2}+32 \delta c_{2}^{2}}{(4c_{2}-\alpha c_{1})^{2}-\delta \alpha^{2}c_{1}^{2}}$\\
$=\frac{\alpha}{4c_{2}-\alpha c_{1}}(1+\frac{32c_{2}^{2}( \frac{ (4c_{2}-\alpha c_{1})^{2}}{8c_{2}(2c_{2}-\alpha c_{1})+(4c_{2}-\alpha c_{1})^{2}})}{(4c_{2}-\alpha c_{1})^{2}- \alpha^{2}c_{1}^{2}(\frac{ (4c_{2}-\alpha c_{1})^{2}}{8c_{2}(2c_{2}-\alpha c_{1})+(4c_{2}-\alpha c_{1})^{2}})})=\frac{\alpha}{4c_{2}-\alpha c_{1}}(\frac{64c_{2}^{2}-16\alpha c_{1}c_{2}}{32c_{2}^{2}-16\alpha c_{1}c_{2}})$\\
$=\frac{\alpha}{2c_{2}-\alpha c_{1}}.$\\
Then above calculations and Theorem 2 imply that  $\overline{x}$ is  the Nash equilibrium level for each player  whenever  $\delta\rightarrow0$, and each player will  choose $\widehat{x}$ for the level of effort if  $\delta\rightarrow \frac{ (4c_{2}-\alpha c_{1})^{2}}{8c_{2}(2c_{2}-\alpha c_{1})+(4c_{2}-\alpha c_{1})^{2}}.$ \\
$\mathbf{Corollary.2}$ In the partnership game, considering the Trigger strategy, for each $\delta \in (0,1)$  the level of effort in an infinitely  repeated game is  determined.




\begin{thebibliography}{99}







\bibitem{Abreu} Abreu, D. (1988 ): On the theory of infinitely repeated games with dis-
counting, Journal of the Econometric Society, JSTOR.


\bibitem{Bernheim} Bernheim, BD. (1984):  Rationalizable strategic behavior, Econometrica: Journal of the Econometric Society,  JSTOR.

\bibitem{Bierman} Bierman, HS.,and  Fernandez, LF. ( 1998 ):  Game theory with economic
applications , Addison -Wesley USA.

\bibitem{Dixit} Dixit, AK. ,and  Skeath, S.  (2015 ): Games of Strategy, Fourth International
Student Edition.

\bibitem{Friedman}Friedman, JW.  (1971): A Non-cooperative Eqililibrium for
Supergames,  Review of Economic Stuadies, 28, 1-12.

\bibitem{Fudenberg}Fudenberg, D., and Maskin, E. (1983):  The Folk Theorem
in Repeated Games with Discounting and Incomplete Information, 
(Working Paper, M.LT.).
\bibitem{Gibbons}Gibbons, R. (1992 ): Game theory for applied economists,Princeton
University Press,Harvester Wheatsheaf.

\bibitem{Green}Green, EJ. , and Porter, RH. (1984): Noncooperative Collusion
under Imperfect Price Information,  Econometrica. 52, 87 100.

\bibitem{GreenE} Green, EJ. (1980): Noncooperative Price-Taking in Large Dynamic
Markets, Journal of Economic Theory, 22,37-64.

\bibitem{Holmstrom} Holmstrom, B. (1982): Moral Hazard in Teams, Bell Journal of
Economics, 13, 324-340.

\bibitem{Lung} Lung, RI., Dumitrescu, D. (2008): Computing Nash Equilibria by Means of
Evolutionary Computation, Communications \&  Control, (suppl. issue)364-8. 

\bibitem{Osborne}  Osborne, MJ. (2004): An introduction to game theory,  Oxford University Press. New York.


\bibitem{OsborneM}  Osborne, MJ. ,and  Rubinstein, A.  (1994): A course in game theory,
MIT Press.

\bibitem{Petrosyan}Petrosyan, LA. ,and  Zenkevich, NA. (1996):  Game Theory (Series on
Optimization, 3), World Scientific Publishers.

\bibitem{Porter} Porter, R H. (1983): Optimal Cartel Trigger-Price Strategies,
Journal of Economic Theory, 29, 313-338.

\bibitem{Radner} Radner R., Myerson R., Maskin E. (1986): An example of a repeated
partnership game with discounting and with uniformly inefficient equilibria, The Review of Economic.

\bibitem{Rubinstein}Rubinstein, A. (1979): Equilibrium in supergames with the overtaking Criterion, Journul of Economic Theory, 21, 1-9.

\bibitem{RubinsteinA} Rubinstein, A.(1977): Equilibrium in Supergames. RM-25. Center for Research in Mathematical Economics and Game Theory, The Hebrew University, Jerusalem. mimeo. 







\end{thebibliography}
\end{document}